\documentclass[reqno]{article}

\usepackage{xcolor}
\usepackage{amssymb}
\usepackage{amsfonts}
\usepackage{amsmath}
\numberwithin{equation}{section}
\usepackage{enumerate}

\newtheorem{theorem}{Theorem}[section]
\newtheorem{proposition}[theorem]{Proposition}
\newtheorem{lemma}[theorem]{Lemma}
\newtheorem{corollary}[theorem]{Corollary}
\newtheorem{definition}[theorem]{Definition}
\newtheorem{example}[theorem]{Example}

\newtheorem{remark}[theorem]{Remark}

\title {Desirable Decompositions of Generalized Nevanlinna Functions}
 \author{Muhamed Borogovac, \\
 muhamed.borogovac@gmail.com}
\begin{document}
\maketitle

\begin{abstract}

For a given generalized Nevanlinna function $Q\in N_{\kappa }\left( H \right)$, we study decompositions that satisfy: $Q=Q_{1}+Q_{2}$; $Q_{i}{\in N}_{\kappa_{i}}\left( H \right)$, and $\kappa_{1}+\kappa_{2}=\kappa $, $0\le \kappa_{i}$, which we call desirable decompositions. In this paper, some sufficient conditions for such decompositions of $Q$ are given. 

One of the main results is a new operator representation of $\hat{Q}\left(z\right):=-{Q(z)}^{-1}$ if $Q\left( z \right):=\Gamma_{0}^{+}\left( A-z\right)^{-1}\Gamma_{0}$, where $A$ is a bounded self-adjoint operator in a Pontryagin space. The new representation is used to get an interesting desirable decomposition of $\hat{Q}$ and to obtain some information about singularities of $\hat{Q}$.

\end{abstract}


\begin{center}
\section{Preliminaries and introduction}\label{s1}
\end{center}
\textbf{1.1.} \textbf{Preliminaries. }Let $N_{0}$, $R$, $C$ denote the sets 
of non-negative integers, real numbers, and complex numbers, respectively. 
Let $(.,.)$ denote a definite scalar product in a Hilbert space $H$ and let 
us denote by $L(H)$ the space of bounded linear operators in $H$.
\\

\begin{definition}\label{def} An operator valued complex function $Q:D\left( Q 
\right)\to L(H)$ belongs to the class of \textbf{generalized Nevanlinna functions} $N_{\kappa }\left( H \right)$ if 
it satisfies the following requirements: 

\begin{itemize}
\item $Q$ is meromorphic in $C\thinspace \backslash \thinspace R$,
\item ${Q\left( \bar{z} \right)}^{\ast }=Q\left( z \right),\thinspace z\in D\left( Q \right),$ 
\item the Nevanlinna kernel $N_{Q}\left( z,\omega \right):=\frac{Q\left( z \right)-{Q\left( \omega \right)}^{\ast }}{z-\bar{\omega }};\thinspace z,\thinspace \omega \in D(Q)\cap \mathbf{C}^{\mathbf{+}}$\textbf{, }
\end{itemize}

has $\kappa $ negative squares, i.e. for arbitrary $n\in N_{0},\thinspace z_{1},\thinspace \mathellipsis ,\thinspace \thinspace z_{n}\in D(Q)\cap \mathbf{C}^{\mathbf{+}}$ and $h_{1},\thinspace \mathellipsis h_{n}\in H$ the Hermitian matrix $\left( N_{Q}\left( z_{i},z_{j} \right)h_{i},h_{j} \right)_{i,j=1}^{n}$ has at most $\kappa $ negative eigenvalues, and for at least one choice of $n;\thinspace z_{1},\thinspace \mathellipsis ,\thinspace \thinspace z_{n}$, and $h_{1},\thinspace \mathellipsis h_{n}$ it has exactly $\kappa $ negative eigenvalues. 
\end{definition}
Let $\kappa \in N_{0}$ and let $\left( K,\thinspace \left[ .,. \right] \right)$ denote a Krein space. If the scalar product $\left[ .,. \right]$ 
has $\kappa (<\infty )$ negative squares it is called a Pontryagin space of index $\kappa $. The definition of a Pontryagin space and other concepts related to it can be found e.g. in \cite{IKL}.

For a bounded linear operator $\Gamma :H\to K$, we denote by $\Gamma ^{+}:K\to H$ the operator defined by $\left( h,\Gamma^{+}k \right):=\left[ \Gamma h,k \right],\thinspace h\in H,\thinspace k\in K$. 

We will deal with the following characterization of generalized Nevanlinna 
functions rather than with Definition \ref{def}.

\begin{proposition}\label{proposition1} A function $\thinspace Q:D(Q)\to L(H)$ is a generalized Nevanlinna function of index $\kappa $, denoted by $Q\in N_{\kappa }(H)$, if and only if it has the representation of the form
\begin{equation}
\label{eq1}
Q\left( z \right)={Q(z_{0})}^{\ast }+(z-\bar{z_{0}})\Gamma^{+}\left( I+\left( z-z_{0} \right)\left( A-z \right)^{-1} \right)\Gamma,\thinspace z\in D\left( Q \right)
\end{equation}
where, A is a self-adjoint linear relation in some Pontryagin space (K, [.,.]) of index $\tilde{\kappa }\ge \kappa ; \Gamma :H\to K$ is a bounded operator; $z_{0}\in \rho \left( A \right)\cap \mathbf{C}^{\mathbf{+}}$ is a fixed point of reference. This representation can be chosen to be minimal, that is
\[
K=cls\left\{ \left( I+\left( z-z_{0} \right)\left( A-z \right)^{-1} \right)\Gamma h:z\in \rho \left( A \right),h\in H \right\}.
\]
The representation is minimal if and only if the negative index of the Pontryagin space $\tilde{\kappa }$ equals $\kappa $. In that case the triplet $(K,\thinspace A,\thinspace \Gamma )$ is uniquely determined (up to isomorphism) and we say that $A$ and $\Gamma $ are closely connected. 
\end{proposition}
Note that in the special case when ``negative index'' $\kappa =0$, the 
Pontryagin space reduces to a Hilbert space. 

Such operator representations were developed by M. G. Krein and H. Langer, 
see e.g. \cite{KL1,KL2} and later translated to representations in terms of 
linear relations (multivalued operators), see e.g. \cite{DLS93}. 

In this paper, a point $\alpha \in C$ is called a generalized pole 
of $Q$ if it is an eigenvalue of the representing relation $A$ of the 
function$\thinspace Q$ given by (\ref{eq1}). It means that it 
may be isolated (i.e. an ordinary pole) as well as an embedded singularity 
of $Q$. 

For later reference we collect the following well known facts into a lemma.

\begin{lemma}\label{lemma1} If $Q\in N_{\kappa }(H)$ is represented by (\ref{eq1}) then it holds 
\begin{equation}
\label{eq2}
{Q\left( z \right)=Q}_{\alpha }(z){+H}_{\alpha }(z),
\end{equation}
where $Q_{\alpha }\in N_{\kappa_{1}}(H), H_{\alpha }\in N_{\kappa_{2}}(H)$ is a holomorphic function at $\alpha , \kappa_{1}+\kappa_{2}=\kappa $. One can always select $Q_{\alpha }$ to be holomorphic at $\infty $. Then $Q_{\alpha}$ admits the representation
\begin{equation}
\label{eq3}
Q_{\alpha }(z)=\Gamma_{0}^{+}\left( A_{0}-z \right)^{-1}\Gamma_{0}\quad ,
\end{equation} 
with a bounded operator $A_{0}$.

If $\Gamma_{0}^{+}\Gamma_{0}$ is not boundedly invertible, one can add a convenient function to $Q_{\alpha }$ in (\ref{eq2}) and subtract the same function from $H_{\alpha }$ so that the new function $Q_{\alpha}^{'}={\Gamma'}_{0}^{+}\left( A_{0}^{'}-z \right)^{-1}{\Gamma'}_{0}$ has the same negative index $\kappa_{1}, {\Gamma'}_{0}^{+}{\Gamma'}_{0}$ is boundedly invertible and $A_{0}^{'}\supseteq A_{0}$. Then the negative index of $H_{\alpha}^{'}$ is $\ge \kappa_{2}.$

By returning to the previous notation, one can consider that $\Gamma_{0}^{+}\Gamma_{0}$ is boundedly invertible in (\ref{eq3}) and now it only holds $\kappa_{1}+\kappa_{2}\ge \kappa $ in (\ref{eq2}). 

In either case, if $\alpha \in R$ is a generalized pole of $Q$, then the operator $A_{0}$ that represents $Q_{\alpha }$ has the same root manifold at $\alpha $ as relation $A$.
\end{lemma}
Because of those properties of $Q_{\alpha }$, it is not a loss of generality 
if one deals with $Q_{\alpha }$ rather than with $Q$ when researching 
properties of $Q$ at $\alpha $. 

Recall here the following statement, see \cite{LaLu,Lu}, which we will also use 
for later references.

\begin{lemma}\label{lemma2} Let the function $Q\in N_{\kappa }(H)$ has minimal representation (\ref{eq1}). If $Q(z_{0})$ is boundedly invertible then the inverse function $\hat{Q}\left( z \right):=-{Q\left( z \right)}^{-1}$ belongs to the class $N_{\kappa }(H)$ and admits the minimal representation 
\begin{equation}
\label{eq4}
\hat{Q}\left( z \right)={-Q\left( \bar{z_{0}} \right)}^{-1}+(z-\bar{z_{0}})\hat{\Gamma }^{+}\left( I+(z-z_{0} \right)\left( \hat{A}-z \right)^{-1})\hat{\Gamma },
\end{equation}
where $\hat{\Gamma }=-\Gamma {Q\left( z_{0} \right)}^{-1}$. Moreover, for $z\in \rho (A)\cap 
\rho (\hat{A})$ and 
\begin{equation}
\label{eq5}
\Gamma_{z}=\left( I+(z-z_{0} \right)\left( A-z \right)^{-1})\Gamma 
\end{equation}
it holds
\begin{equation}
\label{eq6}
\left( \hat{A}-z \right)^{-1}-\left( A-z \right)^{-1}=-\Gamma_{z}{Q\left( z \right)}^{-1}\Gamma_{\bar{z}}^{+}.
\end{equation}
\end{lemma}.
\textbf{1.2. Introduction.} When one studies a complicated object, a way to go is to break it down to simpler components. The same is true with generalized Nevanlinna functions. Various breakdowns of those functions have been proven; some additive (decompositions), some multiplicative (factorizations). In this paper, we will focus on decompositions.

It is well known that a sum $Q$ of generalized Nevanlinna functions that satisfies 

(a) $Q_{i}{\in N}_{\kappa_{i}}(H)\thinspace ,\thinspace 0\le \kappa_{i},\thinspace i=1,\thinspace 2$, 

(b) $Q\left( z \right)=Q_{1}\left( z \right)+Q_{2}\left( z \right)$, 

belongs to some generalized Nevanlinna class $N_{\kappa }\left( H \right)$ and that it holds $\kappa_{1}+\kappa_{2}\ge \kappa $ However, the decompositions with $\kappa_{1}+\kappa_{2}>\kappa $ are not particularly interesting because then the properties of component functions $Q_{i}$ do not add up correctly to the properties of $Q$. In this paper, our main goal is to find necessary and sufficient conditions that functions $Q_{1}$ and $Q_{2}$ satisfying conditions (a) and (b) also satisfy 

(c) $\kappa_{1}+\kappa_{2}=\kappa $. 

A decomposition that satisfy (a), (b) and (c) we call a \textit{desirable  decomposition}. Obviously, that definition can be extended to the sums of more than two functions. 

Some sufficient conditions that a function satisfying (a) and (b) also satisfies (c) were given for scalar functions in \cite{KL2} and for matrix functions in \cite{DaL}, Proposition 3.2. However, those papers only dealt with functions $Q_{i}\in N_{\kappa_{i}}^{nxn}$, $i=1,\thinspace 2$ that have disjoint sets of generalized poles not of positive type, i.e. $\sigma_{0}\left( Q_{1} \right)\cap \sigma_{0}\left( Q_{2} \right)=\emptyset $. In addition, it was assumed: If $\alpha \in \sigma_{0}\left( Q_{j} \right)\cap R$,  then $\lim_{\eta\to 0}{\eta Q_{k}\left( \alpha +i\eta \right)}=0$ and if $\infty \in \sigma _{0}\left( Q_{j} \right)$, then $\thinspace \lim_{\eta\to \infty}{\eta^{-1}Q_{k}\left( i\eta \right)=0}$, $j\ne k,\thinspace j,\thinspace k=1,\thinspace 2$.

In Section \ref{s2}, Theorem \ref{theorem1}, we give some sufficient conditions for desirable decompositions in the most general case, for functions of the form (\ref{eq1}). In addition to that, for a given functions $Q_{i}$ that satisfy (a) we give sufficient conditions that the sum $Q=Q_{1}+Q_{2}$ belongs $N_{\kappa_1+\kappa_{2}}\left(H\right)$, which means that the number of negative squares is preserved. In order to do that we had to introduce the following assumption

(d) $\Gamma^{+}$ is injection.

Our results also apply to desirable decompositions of a given function $Q$ where components $Q_i$ have the same critical point. A decomposition where the decomposing functions have a common critical point we call a decomposition \textit{within a critical point}.

Example \ref{example1} is complementary to the statements \ref{lemma3} through \ref{proposition2} because it explains assumptions of those statements. It also shows us that converse statement of Theorem \ref{theorem1} (i) does not hold. 

In the short Section \ref{s3}, we derive one desirable decomposition within a 
generalized pole $\alpha $ in terms of the maximal Jordan chains in that pole. 

In Section \ref{s4}, the main result is Theorem \ref{theorem2}. In that theorem we assume $Q\left( z \right):=\Gamma_{0}^{+}\left( A-z \right)^{-1}\Gamma_{0}\in N_{\kappa }(H)$, where $A$ is a bounded self-adjoint operator in a Pontryagin space and $\Gamma_{0}^{+}\Gamma_{0}$ is boundedly invertible and we prove a new, operator representation of $\hat{Q}\left( z \right):=-{Q(z)}^{-1}$. 

In Section \ref{s5}, we use that representation of $\hat{Q}$ to find the 
decomposition 
\[
\hat{Q}(z)=\hat{Q}_{1}(z)+\hat{Q}_{2}(z);
\quad
\hat{Q}_{1}\in N_{\hat{\kappa }_{1}}(H),
\quad
\hat{Q}_{2}\in N_{\hat{\kappa }_{2}}\left( H \right);
\quad
\hat{\kappa }_{1}+\hat{\kappa }_{2}=\kappa ,
\]
where one of the components, e.g. $\hat{Q}_{1}$ can have only a zero in the 
critical point $\alpha $ of $\hat{Q}$. 

\section{Desirable decomposition of a generalized Nevanlinna function}\label{s2}
\textbf{2.1.} According to Lemma \ref{lemma1}, representations of the form (\ref{eq3}) play an 
important role in decompositions of generalized Nevanlinna functions. To simplify notation, we will deal with a function $Q\in N_{\kappa}(H)$ 
\begin{equation}
\label{eq7}
Q(z)=\Gamma_{0}^{+}\left( A-z \right)^{-1}\Gamma_{0}\mathrm{\thinspace },
\end{equation}
where $\Gamma_{0}:H\to K$ is a bounded operator and $A$ is a bounded self-adjoint operator in a Pontryagin space $K$. We will always denote by $\Gamma $ operator used in representation (\ref{eq1}) and by $\Gamma_{0}$ operator used in the special case, 
representation (\ref{eq7}).
\\

\begin{lemma}\label{lemma3} 

\begin{enumerate}[(i)]%, (i), (ii), ...
\item If $\Gamma^{+}$ in representation (\ref{eq1}) is an injection, then (\ref{eq1}) is minimal representation of $Q$. 
\item If $\Gamma_{0}^{+}\Gamma_{0}$ is injection, then the function $Q$ given by (\ref{eq7}) satisfies
\begin{equation}
\label{eq8}
\left( f,Q\left( z \right)h \right)=0,\thinspace \forall z\in D\left( Q\right),\thinspace \forall h\in H\to f=0
\end{equation}
\item If (\ref{eq8}) holds then $\Gamma_{0}$ is injection.
\item Assume (\ref{eq7}) is minimal. If $\Gamma_{0}$ is injection then (\ref{eq8}) holds. 
\end{enumerate}
\end{lemma}
\textbf{Proof.} (i). Assume $\Gamma^{+}$ is an injection and 
\[
0=\left[ y,\left( I+\left( z-z_{0} \right)\left( A-z \right)^{-1}\right)\Gamma h \right],\thinspace \forall z\in \rho \left( A \right), \forall h\in H
\]
Then 
\[
\left( \Gamma^{+}\left( I+\left( \bar{z}-\bar{z_{0}} \right)\left( 
A-\bar{z} \right)^{-1} \right)y,h \right)=0,\forall z\in \rho \left( A 
\right),\thinspace \forall h\in H
\]
As $\Gamma^{+}$ is injection, we conclude
\begin{equation}
\label{eq9}
\left( I+\left( \bar{z}-\bar{z_{0}} \right)\left( A-\bar{z} \right)^{-1} 
\right)y=0
\end{equation}
The obvious solution of the equation (\ref{eq9}) is $y=0$. If that is the only 
solution of (\ref{eq9}) then minimality of the representation (\ref{eq1}) is proven. 
Assume to the contrary that equation (\ref{eq9}) has a nonzero solution $y$. Then it follows: 

As $\bar{z}\in \rho \left( A \right)$, the relation $\quad \left( A-\bar{z} \right)^{-1}$ is 
defined on $K$. That implies, 
\[
\left( A-\bar{z} \right)\left( A-\bar{z} \right)^{-1}\supseteq I
\]
Then from (\ref{eq9}) we have $-\left( \bar{z}-\bar{z_{0}} \right)y\in \left( 
A-\bar{z} \right)y$ and therefore $\bar{z_{0}}y\in Ay$. This means that $\bar{z_{0}}$ is an eigenvalue of $A$. It is in contradiction with the fact that $\bar{z_{0}}$ is regular point of $A$ as the symmetrical point of $z_{0}\in C^{+}\cap \rho \left( A \right)$. Therefore, it has to be $y=0$, which proves that representation (\ref{eq1}) is minimal. 

(ii) $0=\left( f,Q\left( z \right)h \right),\thinspace \forall z\in 
D\left( Q \right),\thinspace \forall h\in H$
\[
\to 0=\left( f,zQ\left( z \right)h \right),\thinspace \forall z\in 
D\left( Q \right),\thinspace \forall h\in H
\]
\[
\to 0=\lim_{z\to \infty}{\left( f,zQ\left( z \right)h \right)=-\left( f,\Gamma_{0}^{+}\Gamma 
_{0}h \right)=-\left( \Gamma_{0}^{+}\Gamma_{0}f,h \right),\thinspace 
\forall h\in H.}
\]
As $H$ is a Hilbert space and $\Gamma_{0}^{+}\Gamma_{0}$ is injection, it 
follows $f=0$.

(iii) Assume that (\ref{eq8}) holds and $\Gamma_{0}f=0$ . Then,
\[
0=\left[ \Gamma_{0}f,\Gamma_{0}h \right]=\left[ \Gamma_{0}f,\left( A-z 
\right)^{-1}\Gamma_{0}h \right]=\left( f,Q\left( z \right)h \right),\thinspace \forall z\in D\left( Q \right),\thinspace \forall h\in H.
\]
Now from (\ref{eq8}) it follows $f=0$, which proves that $\Gamma_{0}$ is 
injection. 

By similar arguments the statement (iv) can be proven.

Note that the representation (\ref{eq7}) can always be selected to be minimal. 
Then, the statements (iii) and (iv) mean that (\ref{eq8}) holds if and only if 
$\Gamma_{0}$ is injection. 

The converse statement of (iv) does not hold. In the Example \ref{example1} we will see 
that it is possible to have (\ref{eq8}) satisfied and $\Gamma_{0}$ is injection, but the corresponding representation (\ref{eq7}) does not need to be minimal. 
\\

\textbf{2.2.} Let us now assume that functions $Q_{i}\in N_{\kappa 
_{i}}(H)$ are of the form (\ref{eq1}) represented by triplets $\left( 
K_{i},A_{i},\Gamma_{i} \right),\thinspace \thinspace i=1,\thinspace 2$, 
i.e. 
\begin{equation}
\label{eq10}
Q_{i}\left( z \right)={Q_{i}(z_{0})}^{\ast }+(z-\bar{z_{0}})\Gamma 
_{i}^{+}\left( I_{i}+\left( z-z_{0} \right)\left( A_{i}-z \right)^{-1} 
\right)\Gamma_{i}
\end{equation}
where $A_{i}$ are self-adjoint relations and $\Gamma_{i}:H\to K_{i}$ . Then we can 
introduce the triplet $(\tilde{K},\thinspace A,\thinspace \Gamma )$ by:
\\
\begin{center}
$\tilde{K}:=K_{1}\left[ + \right]K_{2}$,
\end{center}
\begin{equation}
\label{eq11}
A:=\left\{ \left\{ k_{1}\left[ + \right]k_{2},h_{1}\left[ + \right]h_{2} 
\right\}:\left\{ k_{i},h_{i} \right\}\in A_{i},\thinspace \thinspace 
i=1,\thinspace 2 \right\}
\end{equation} 
Because $A_{i}$ are self-adjoint in $K_{i}$, $A$ is self-adjoint in 
$K_{1}\left[ + \right]K_{2}$, i.e. $A=A^{[\ast ]}$.
\begin{equation}
\label{eq12}
\Gamma :H\to K_{1}\left[ + \right]K_{2},\thinspace \Gamma \left( h 
\right):=\Gamma_{1}\left( h \right)\left[ + \right]\Gamma_{2}\left( h 
\right),\thinspace \Gamma_{i}\left( h \right)\in K_{i}\thinspace 
\end{equation}
Then ${\Gamma^{+}=\Gamma }_{1}^{+}+\Gamma_{2}^{+}:K_{1}\left[ + 
\right]K_{2}\to H$ is defined by 
\begin{equation}
\label{eq13}
\left[ \Gamma h,k_{1}\left[ + \right]k_{2} \right]=(h,\Gamma 
_{1}^{+}k_{1}+\Gamma_{2}^{+}k_{2})
\end{equation}
where $k_{i}\in K_{i}$.

The same definitions hold when functions $Q_{i}$ are of the form (\ref{eq7}) with 
$\Gamma_{0i}$ and $\Gamma_{0}$ replacing $\Gamma_{i}$ and $\Gamma $. The above definitions prepared us for the following lemma. 

\begin{lemma}\label{lemma4} Let functions $Q_{i}$ be represented by (\ref{eq10}). For the 
function $Q:=Q_{1}+Q_{2}$, let us observe the representation
\[
Q\left( z \right)={Q_{1}\left( z_{0} \right)}^{\ast }+{Q_{2}\left( z_{0} 
\right)}^{\ast }+
\]
\begin{equation}
\label{eq14}
+(z-\bar{z_{0}})\Gamma^{+}\left( {\begin{array}{*{20}c}
{I_{1}+(z-z_{0})\left( A_{1}-z \right)}^{-1} & 0\\
0 & {I_{2}+(z-z_{0})\left( A_{2}-z \right)}^{-1}\\
\end{array} } \right)\Gamma 
\end{equation}
where functions $\Gamma $ and $\Gamma^{+}$ are defined by (\ref{eq12}), (\ref{eq13}), respectively. 

\begin{enumerate}[(i)]%, (i), (ii), ...
\item If operator $\Gamma^{+}:K_{1}\left[ + \right]K_{2}\to H$ introduced by (\ref{eq13}) is injection, then (\ref{eq14}) is minimal.
\item Let now functions $Q_{i},\thinspace i=1,\thinspace 2$ be given by (\ref{eq7}), let $\Gamma_{0}$ and $\Gamma_{0}^{+}$ be defined by (\ref{eq12}), (\ref{eq13}) and let representation 
\end{enumerate}
\begin{equation}
\label{eq15}
Q(z):=Q_{1}(z)+Q_{2}(z)=\Gamma_{0}^{+}\left( {\begin{array}{*{20}c}
\left( A_{1}-z \right)^{-1} & 0\\
0 & \left( A_{2}-z \right)^{-1}\\
\end{array} } \right)\Gamma_{0}
\end{equation}
be minimal. If at least one of $\Gamma_{0i}$ is an injection then it holds
\[
\left( f,Q\left( z \right)h \right)=0,\thinspace \thinspace \forall z\in 
D\left( Q \right),\thinspace \forall h\in H\to f=0
\]
\end{lemma}
\textbf{Proof.} (i).The statement follows when we apply Lemma \ref{lemma3} (i) to 
function $Q$ represented by (\ref{eq14}), where $K$, $A$, $\Gamma $ and $\Gamma 
^{+}$are defined by (\ref{eq11}), (\ref{eq12}) and (\ref{eq13}). 

(ii) Assume now 
\[
\left( f,Q\left( z \right)h \right)=\left( f,Q_{1}\left( z \right)h 
\right)+\left( f,Q_{2}\left( z \right)h \right)=0,\forall z\in D\left( Q 
\right),\thinspace \forall h\in H
\]
Then
\[
\left[ \Gamma_{01}f,\left( A_{1}-z \right)^{-1}\Gamma_{01}h \right]+\left[ 
\Gamma_{02}f,\left( A_{2}-z \right)^{-1}\Gamma_{02}h \right]=0,\thinspace 
\forall z\in \rho \left( A \right),\thinspace \forall h\in H
\]
From assumption of minimality of $Q$ we conclude
\[
\left( {\begin{array}{*{20}c}
\Gamma_{01}f\\
\Gamma_{02}f\\
\end{array} } \right)=0
\]
From the assumption that at least one $\Gamma_{0i}$ is injection we have 
$f=0$. 
\\

Note, we did \textbf{not} assume minimality of individual representations (\ref{eq10}). Minimality of components does not guarantee the minimality of the sum (\ref{eq14}), while the minimality of the sum guaranties minimality of the components in the representations (\ref{eq14}) and (\ref{eq15}).

\begin{theorem}\label{theorem1} (i) Assume representation (\ref{eq1}) of $Q\in N_{\kappa }(H)$ is minimal, and there exist non-degenerate, invariant with respect to $A$\, sub-spaces $K_{i}, K_{1}\left[ + \right]K_{2}=K$. Then
\begin{enumerate}[(a)]%, (a), (b), ...
\item $\exists Q_{i}{\in N}_{\kappa_{i}}\left( H \right),$ minimally represented by triplets $\left( K_{i},A_{i},\Gamma_{i}\thinspace \right), $
\item $Q\left( z \right)=Q_{1}\left( z \right)+Q_{2}\left( z \right),$
\item $\thinspace \kappa_{1}+\kappa_{2}=\kappa .$ 

(ii) If conditions (a), (b) and 
\item ${\Gamma^{+}=\Gamma _{1}^{+}+\Gamma_{2}^{+}}$ is an injection;  
\end{enumerate}
are satisfied, then the unique minimal representation of $Q$ is given by (\ref{eq14}), where the representing triplet $(K,\thinspace A,\thinspace \Gamma )$ is defined by (\ref{eq11}), (\ref{eq12}) and (\ref{eq13}). In addition, $Q\in N_{\kappa_{1}+\kappa_{2}}(H)$, i.e. (c) holds.

\end{theorem}
\textbf{Proof (i)} We will prove the proposition under seemingly more general assumptions. We will assume existence of only one non-degenerate invariant subspace $K_{1}\subseteq K$. Then we introduce the orthogonal projection onto \quad $K_{1}$, \quad $E_{1}:K\to K_{1}.$ For \quad $E_{2}:=I-E_{1}$ and \quad $K_{2}:=E_{2}K$ 
the following decomposition holds
\\

 $K=K_{1}\left[ + \right]K_{2}$,
\\

where $K_{1}$ and $K_{2}$ are Pontryagin subspaces of negative indexes 
$\kappa_{i}$, $0\le \kappa_{i},\thinspace \thinspace i=1,2$, $\kappa 
_{1}+\kappa_{2}=\kappa $. 

Because $A$ is a self-adjoint relation, $K_{2}$ is also invariant with 
respect to $A$. Then 
\[
I+\left( z-z_{0} \right)\left( A-z \right)^{-1}=
\]
\begin{equation}
\label{eq16}
=\left( {\begin{array}{*{20}c}
E_{1}+{\left( z-z_{0} \right)[E_{1}\left( A-zI \right)E_{1}]}^{-1} & 0\\
0 & E_{2}+{\left( z-z_{0} \right)[E_{2}\left( A-zI \right)E_{2}]}^{-1}\\
\end{array} } \right).
\end{equation}
If we introduce $A_{i}:=E_{i}AE_{i}:K_{i}\to K_{i}$, $\Gamma 
_{i}:=E_{i}\Gamma $, $i=1,\thinspace 2$, then it holds
\[
Q\left( z \right)=Q_{1}\left( z \right)+Q_{2}\left( z \right),
\]
where
\begin{equation}
\label{eq17}
Q_{i}\left( z \right):={Q_{i}\left( z_{0} \right)}^{\ast }+\left( 
z-\bar{z_{0}} \right)\Gamma_{i}^{+}\left( {I+\left( z-z_{0} \right)\left( 
A_{i}-z \right)}^{-1} \right)\Gamma_{i}.
\end{equation}
From (\ref{eq16}) and from the minimality of the representation (\ref{eq1}), the minimality of representations (\ref{eq17}) follows. 

Indeed, for $y_{1}[+]y_{2}\in K_{1}[+]K_{2}$ minimality of (\ref{eq1}) means 
\[
\left[ \left( {\begin{array}{*{20}c}
y_{1}\\
y_{2}\\
\end{array} } \right),\left( {\begin{array}{*{20}c}
E_{1}+{\left( z-z_{0} \right)[E_{1}\left( A-zI \right)E_{1}]}^{-1} & 0\\
0 & E_{2}+{\left( z-z_{0} \right)[E_{2}\left( A-zI \right)E_{2}]}^{-1}\\
\end{array} } \right)\left( {\begin{array}{*{20}c}
\Gamma_{1}h\\
\Gamma_{2}h\\
\end{array} } \right) \right]
\]
\[
\thinspace =0, \forall z\in \rho \left( A \right),\thinspace \forall 
h\in H,\to \left( {\begin{array}{*{20}c}
y_{1}\\
y_{2}\\
\end{array} } \right)=0
\]
If we keep $y_{2}=0$, we conclude that $Q_{1}$ is minimally represented by 
$\left( K_{1},A_{1},\Gamma_{1}\thinspace \right)$. By the same token we 
conclude that $Q_{2}$ is minimally represented by $\left( K_{2},A_{2},\Gamma 
_{2}\thinspace \right)$.

It further means that negative indexes of $Q_{i}$ and $K_{i}$ are equal. 
Hence, $Q_{i}{\in N}_{\kappa_{i}}(H),\thinspace i=1,2$ and from $K=K_{1}\left[ + \right]K_{2}$ we get $\kappa_{1}+\kappa_{2}=\kappa $. 

That proves (a), (b) and (c). 

\textbf{Proof (ii)} Assume now that functions \quad $Q_{i}$, represented by (\ref{eq1}), satisfy conditions (a), (b) and (d). According to Lemma \ref{lemma4} (i), the representation (\ref{eq14}) is minimal representation in terms of the triplet 
$(\tilde{K},\thinspace A,\thinspace \Gamma )$ defined by (\ref{eq11}), (\ref{eq12}) and 
(\ref{eq13}). The minimality of representations of $Q_{i}$ implies that negative 
indexes of $K_{i}$ are $\kappa_{i}$, respectively. From $\tilde{K}=K_{1}\left[ + \right]K_{2}$ it follows that negative index $\tilde{\kappa }$ of $\tilde{K}$ 
satisfies $\tilde{\kappa }=\kappa_{1}+\kappa_{2}$.

From the minimality of the representation (\ref{eq14}) of $Q$ in terms of $(\tilde{K},\thinspace A,\thinspace \Gamma )$ it follows $Q\in N_{\tilde{\kappa }}\left( H \right)$. Hence, $Q\in N_{\kappa_{1}+\kappa_{2}}\left( H \right)$.

Then $K_{1}\subseteq \tilde{K}$ is the non-degenerate subspace invariant with respect to relation $A$. 

From the uniqueness (up to isomorphism) of the representing triplet of the minimal representation, we conclude that representation (\ref{eq1}) is of the form (\ref{eq14}), and we can denote $\tilde{K}$ by $K$. This proves the statement (ii).
\\

By means of the triplet $(\tilde{K},\thinspace A,\thinspace \tilde{\Gamma 
})$, where $\tilde{K}$ and $A$ are as before and $\tilde{\Gamma 
}:H_{1}\left( + \right)H_{2}\to K_{1}\left[ + \right]K_{2}$ is defined by 
\[
\tilde{\Gamma }\left( h_{1}(+)h_{2} \right):=\Gamma_{1}h_{1}[+]\Gamma 
_{2}h_{2}
\]
it is easy to prove the following proposition. 

\begin{proposition}\label{proposition2} If $Q_{i}{\epsilon N}_{\kappa_{i}}\left( H_{i} 
\right),\thinspace i=1,\thinspace 2$\textit, and
\[
\tilde{Q}(z):=\left( {\begin{array}{*{20}c}
Q_{1}\left( z \right) & 0\\
0 & Q_{2}\left( z \right)\\
\end{array} } \right),
\]
then $\tilde{Q}\in N_{\kappa_{1}+\kappa_{2}}\left( H_{1}\left( + 
\right)H_{2} \right). $
\end{proposition}

\textbf{2.3.} The following example explains many assumptions of the statements \ref{lemma3} through \ref{proposition2}, making them natural. 

\begin{example}\label{example1}Given the matrix function
\[
Q\left( z \right)=-\left( {\begin{array}{*{20}c}
z^{-1}+z^{-2} & z^{-1}\\
z^{-1} & z^{-1}\\
\end{array} } \right).
\]
\end{example}
Obviously $Q$ is a regular (boundedly invertible for every $z\ne 0)$ 
function of the form (\ref{eq7}) but $Q^{'}\left( \infty \right)=-\Gamma 
_{0}^{+}\Gamma_{0}$ is not even an injection. 
\[
Q\left( z \right)=-\left( {\begin{array}{*{20}c}
z^{-1}+z^{-2} & z^{-1}\\
z^{-1} & z^{-1}\\
\end{array} } \right)
\]
\[
=-\left( {\begin{array}{*{20}c}
z^{-1} & 0\\
0 & z^{-1}\\
\end{array} } \right)-\left( {\begin{array}{*{20}c}
z^{-2} & z^{-1}\\
z^{-1} & 0\\
\end{array} } \right)=:Q_{1}(z)+Q_{2}(z).
\]
Because the scalar function $q\left( z \right)=-z^{-1}$ belongs to class 
$N_{0}$, according to Proposition \ref{proposition2} it holds 
\[
Q_{1}\left( z \right):=-\left( {\begin{array}{*{20}c}
z^{-1} & 0\\
0 & z^{-1}\\
\end{array} } \right)=\Gamma_{1}^{+}\left( A_{1}-z \right)^{-1}\Gamma 
_{1}\in N_{0}^{2x2},
\]
where $\Gamma_{1}=\Gamma_{1}^{+}=\Gamma_{1}^{\ast }=\left( 
{\begin{array}{*{20}c}
1 & 0\\
0 & 1\\
\end{array} } \right)$, $A_{1}=\left( {\begin{array}{*{20}c}
0 & 0\\
0 & 0\\
\end{array} } \right)$.
\[
Q_{2}\left( z \right):=-\left( {\begin{array}{*{20}c}
z^{-2} & z^{-1}\\
z^{-1} & 0\\
\end{array} } \right)=\Gamma_{2}^{+}\left( A_{2}-z \right)^{-1}\Gamma 
_{2},
\]
where $A_{2}=\left( {\begin{array}{*{20}c}
0 & 1\\
0 & 0\\
\end{array} } \right)$, $J_{2}=\left( {\begin{array}{*{20}c}
0 & 1\\
1 & 0\\
\end{array} } \right)$, $\Gamma_{2}=\left( {\begin{array}{*{20}c}
0 & 1\\
1 & 0\\
\end{array} } \right)$.

It is easy to see that both functions $Q_{i}$ are minimally represented. 

Formula (\ref{eq15}) gives here
\[
Q_{1}\left( z \right)+Q_{2}\left( z \right)=\left( \Gamma_{1}^{+}\thinspace 
\thinspace \thinspace \thinspace \Gamma_{2}^{+} \right)\left( 
{\begin{array}{*{20}c}
\left( A_{1}-z \right)^{-1} & 0\\
0 & \left( A_{2}-z \right)^{-1}\\
\end{array} } \right)\left( \frac{\Gamma_{1}}{\Gamma_{2}} \right)=
\]
\begin{equation}
\label{eq18}
=-\left( {\begin{array}{*{20}c}
\left( {\begin{array}{*{20}c}
1 & 0\\
0 & 1\\
\end{array} } \right) & \left( {\begin{array}{*{20}c}
1 & 0\\
0 & 1\\
\end{array} } \right)\\
\end{array} } \right)\thinspace \left( {\begin{array}{*{20}c}
\left( {\begin{array}{*{20}c}
z^{-1} & 0\\
0 & z^{-1}\\
\end{array} } \right) & \left( {\begin{array}{*{20}c}
0 & 0\\
0 & 0\\
\end{array} } \right)\\
\left( {\begin{array}{*{20}c}
0 & 0\\
0 & 0\\
\end{array} } \right) & \left( {\begin{array}{*{20}c}
z^{-1} & z^{-2}\\
0 & z^{-1}\\
\end{array} } \right)\\
\end{array} } \right)\left( \frac{\left( {\begin{array}{*{20}c}
1 & 0\\
0 & 1\\
\end{array} } \right)}{\left( {\begin{array}{*{20}c}
0 & 1\\
1 & 0\\
\end{array} } \right)} \right)
\end{equation}
\[
=-\left( {\begin{array}{*{20}c}
z^{-1}+z^{-2} & z^{-1}\\
z^{-1} & z^{-1}\\
\end{array} } \right).
\]
Let $(\tilde{K},\thinspace A,\thinspace \Gamma_{0})$ be the triplet created 
by means formulae (\ref{eq11}), (\ref{eq12}) and (\ref{eq13}).
\[
J\mathrm{=}\left( {\begin{array}{l}
 {\begin{array}{*{20}c}
\mathrm{1} & \mathrm{0}\\
\mathrm{0} & \mathrm{1}\\
\mathrm{0} & \mathrm{0}\\
\end{array} }\mathrm{\thinspace \thinspace \thinspace 
}{\begin{array}{*{20}c}
\mathrm{0} & \mathrm{0}\\
\mathrm{0} & \mathrm{0}\\
\mathrm{0} & \mathrm{1}\\
\end{array} } \\ 
 {\begin{array}{*{20}c}
\mathrm{0} & \mathrm{0}\\
\end{array} }\mathrm{\thinspace \thinspace }{\begin{array}{*{20}c}
\mathrm{1} & \mathrm{0}\\
\end{array} } \\ 
 \end{array}} \right),
\quad
A\mathrm{:=}\left( {\begin{array}{l}
 {\begin{array}{*{20}c}
\mathrm{0} & \mathrm{0}\\
\mathrm{0} & \mathrm{0}\\
\mathrm{0} & \mathrm{0}\\
\end{array} }\mathrm{\thinspace \thinspace \thinspace 
}{\begin{array}{*{20}c}
\mathrm{0} & \mathrm{0}\\
\mathrm{0} & \mathrm{0}\\
\mathrm{0} & \mathrm{1}\\
\end{array} } \\ 
 {\begin{array}{*{20}c}
\mathrm{0} & \mathrm{0}\\
\end{array} }\mathrm{\thinspace \thinspace }{\begin{array}{*{20}c}
\mathrm{0} & \mathrm{0}\\
\end{array} } \\ 
 \end{array}} \right),
\]
\[
\Gamma_{0}\mathrm{=}\left( {\begin{array}{l}
 {\begin{array}{*{20}c}
\mathrm{1} & \mathrm{0}\\
\mathrm{0} & 1\\
\mathrm{0} & \mathrm{1}\\
\end{array} } \\ 
 {\begin{array}{*{20}c}
\mathrm{1} & \mathrm{0}\\
\end{array} } \\ 
 \end{array}} \right),
\quad
\Gamma_{0}^{+}\mathrm{=}\left( {\begin{array}{*{20}c}
\mathrm{1} & 0 & \mathrm{1\thinspace }\\
\mathrm{0} & 1 & \mathrm{0}\\
\end{array} }\mathrm{\thinspace \thinspace \thinspace 
}{\begin{array}{*{20}c}
\mathrm{0}\\
\mathrm{1}\\
\end{array} } \right).
\]
Obviously, $\Gamma_{0}$ is injection but $\Gamma_{0}^{+}$ and $\Gamma 
_{0}^{+}\Gamma_{0}$ are not. Assume now
\[
\left( \left( {\begin{array}{*{20}c}
f_{1}\\
f_{2}\\
\end{array} } \right),\left( {\begin{array}{*{20}c}
z^{-1}+z^{-2} & z^{-1}\\
z^{-1} & z^{-1}\\
\end{array} } \right)\left( {\begin{array}{*{20}c}
h_{1}\\
h_{2}\\
\end{array} } \right) \right)=0,\thinspace \forall h=\left( 
{\begin{array}{*{20}c}
h_{1}\\
h_{2}\\
\end{array} } \right)\in \mathbf{C}^{2}
\]
As the functions $g_{1}(z):=\left( z^{-1}+z^{-2} \right)h_{1}+z^{-1}h_{2}$ 
and $g_{2}(z):=z^{-1}h_{1}+z^{-1}h_{2}$ are obviously linearly independent 
for every fixed $h=\left( {\begin{array}{*{20}c}
h_{1}\\
h_{2}\\
\end{array} } \right)$, it has to be $\left( {\begin{array}{*{20}c}
f_{1}\\
f_{2}\\
\end{array} } \right)=0$. This means that (\ref{eq8}) holds and $\Gamma_{0}$ is injection. However, it is easy to verify that the corresponding representation (\ref{eq18}) of $Q$ is not minimal. Hence, converse statements of Lemma \ref{lemma3} (iv) and Lemma \ref{lemma4} (ii) do not hold. 

Note that conditions (a), (b) and (c) are satisfied in this example but it is not sufficient for minimality of (\ref{eq18}). That means that the converse statement of Theorem \ref{theorem1} (i) does not hold. It justifies introduction of the condition (d) in the study of desirable decompositions. Note also that minimal representation of $Q$ must be different form (\ref{eq18}).
\section{Decomposition of the Pontryagin space by means of Jordan chains of a self-adjoint relation}\label{s3}
\textbf{3.1.} Let us denote root manifold (algebraic eigenspace) of the representing relation $A$ at $\alpha \in R$ by $S_{\alpha }(A):=\left\{x:\exists r\in N,\thinspace \thinspace \left( A-\alpha \right)^{r}x=0 \right\}$. 

Let $X=\left\{ x_{k}\thinspace ,\thinspace k=0,\thinspace \mathellipsis ,\thinspace l-1 \right\}$ be a maximal non-degenerate Jordan chain at $\alpha$ of the representing relation $A$ of $Q\in N_{\kappa }(H)$. 
According to Lemma \ref{lemma1}, $X$ is also Jordan chain of the bounded self-adjoint 
operator $A_{0}$ representing the function $Q_{\alpha }$. 

\begin{proposition}\label{proposition3} Let $Q\in N_{\kappa }(H)$ be given by minimal representation (\ref{eq1}) and let $X$ be a maximal non-degenerate Jordan chain of the length $l$ of the representing relation A at $\alpha$. Then 

\begin{enumerate} [(i)]%, (i), (ii), ...
\item $Q\left( z \right)=Q_{1}\left( z \right)+Q_{2}\left( z \right),$ 

where $Q_{i}{\in N}_{\kappa_{i}}\left( H \right), i=1,\thinspace 2; \thinspace 
\kappa_{1}+\kappa_{2}=\kappa ; \sigma \left( Q_{1} \right)=\left\{ \alpha 
\right\}$ and $Q_{1}\left( z \right)=\Gamma_{1}^{+}\left( A_{1}-z 
\right)^{-1}\Gamma_{1} . $
\item There exist $x_{l-1}\in \Gamma_{1}(H)$, such that
\end{enumerate}
\[
X=\left\{ \left( A-\alpha \right)^{l-1-k}x_{l-1}\thinspace ,\thinspace 
k=0,\thinspace \mathellipsis ,\thinspace l-1 \right\}.
\]
\end{proposition}
\textbf{Proof.} \textbf{(i)} Let us introduce $S_{\alpha }\left( x_{0} 
\right):=c.l.s.(X)$ and let $E:K\to S_{\alpha }(x_{0})$ denote the orthogonal 
projection onto \quad $S_{\alpha }(x_{0})$. We can define $K_{1}:=\thinspace 
S_{\alpha }(x_{0})$. Then the statement (i) follows directly from Lemma \ref{lemma1} and Theorem \ref{theorem1}. 

\textbf{(ii)} As before$\thinspace \Gamma_{1}:=E\Gamma$ and $\mathrm{\thinspace 
}A_{1}:=AE=EAE$ are closely connected in $\thinspace E(K)$. As $A_{1}$ is 
bounded operator, and $A_{1}$ and $\Gamma_{1}$ are closely connected, we 
have 
\[
S_{\alpha }(x_{0})=c.l.s.\left\{ A^{i}\Gamma_{1}\left( H \right),\thinspace 
i=0,\thinspace 1,\thinspace \mathellipsis \right\}=c.l.s.\left\{ {(A-\alpha 
I)}^{i}\Gamma_{1}\left( H \right),\thinspace i=0,\thinspace 1,\thinspace 
\mathellipsis \right\}.
\]
Therefore, the last vector in the given Jordan chain, $x_{l-1}$ must have a 
representation of the form 
\[
x_{l-1}=\Gamma_{1}h_{0}+(A-\alpha I)y\quad ,
\]
where $y:=\sum\limits_{i=0}^\infty {{(A-\alpha I)}^{i}\Gamma_{1}h_{i+1}} 
\in S_{\alpha }(x_{0})$.

Obviously, ${(A-\alpha I)}^{i}\Gamma_{1}h_{i+1}=0$ for every $i\ge l$.

If $y=0$, then $x_{l-1}=\Gamma_{1}h_{0}$, which proves (ii). If $y\ne 0$, 
then it follows 
\[
x_{0}=\left( A-\alpha I \right)^{l-1}x_{l-1}=\left( A-\alpha I 
\right)^{l-1}\Gamma_{1}h_{0}.
\]
Hence, we can take $x_{l-1}=\Gamma_{1}h_{0}\in E\Gamma (H)$.

\begin{remark}\label{remark2} Obviously, a typical situation is that $\Gamma 
_{1}\left( H \right)$ is a proper subset of $K_{1}$. There exist examples of 
maximal Jordan chains that do not begin at $\Gamma_{1}\left( H \right)$ 
(meaning $x_{l-1}\notin \Gamma_{1}(H))$. Also, there are examples of chains 
that begin in $\Gamma_{1}(H)$ that are not maximal. The meaning of the 
statement (ii) is that structure of the algebraic eigen-space $S_{\alpha}(A)$ can be characterized by means of maximal non-degenerate Jordan chain with the last vectors $x_{l-1}\in \Gamma_{1}\left( H \right)$. 
\end{remark}

\textbf{3.2.} Let $\alpha \in R$ be a generalized pole not of positive 
type of $Q\in N_{\kappa }(H)$. We will focus on the decomposition 
within a single critical point. Therefore, it is not a loss of generality to 
assume that $Q$ admits representation (\ref{eq3}) and that $\alpha \in R$ is 
the single critical point. For simplicity, we again use $Q$, $\thinspace A$ 
and $\Gamma $ rather than of $Q_{\alpha }$, $A_{0}$ and $\Gamma_{0}$. 

Let $K_{0}\subseteq K$ be the Hilbert subspace that consists of all positive 
eigenvectors of the representing operator $A$ et $\alpha $. Obviously, 
$K_{0}$ is invariant subspace with respect to $A$. Let $E_{0}:K\to K_{0}$ be 
the orthogonal projection and $\Gamma_{0}:=E_{0}\Gamma $. Then the 
Pontryagin space $\left( I-E_{0} \right)K$ is also invariant with respect to 
$A$, and operators $\tilde{\Gamma }:=\left( I-E_{0} \right)\Gamma $ and 
$\tilde{A}:=\left( I-E_{0} \right)A\left( I-E_{0} \right)$ are closely 
connected. 

Now let $x_{0}^{1},\thinspace \mathellipsis ,\thinspace x_{l_{1}-1}^{1}$ be a 
maximal non-degenerate Jordan chain of $\tilde{A}$ at 
$\alpha $ in the Pontryagin space $\left( I-E_{0} \right)K$. We define: $E_{1}:\left( I-E_{0} \right)K\to S_{\alpha }(x_{0}^{1})$; 
${K_{1}:=E}_{1}\left( I-E_{0} \right)K$. Then $A_{1}={AE}_{1}$ and \quad $\Gamma 
_{1}:=E_{1}\Gamma $ are closely connected,$\thinspace \kappa_{1}$ is index 
of the Pontryagin space $K_{1}$. According to Proposition \ref{proposition3} (ii) we can 
consider $x_{l_{1}-1}^{1}=\Gamma_{1}h_{1}$. We continue that process until 
we exhaust all non-degenerate Jordan chains. Assume that there are 
$r>0$ such chains. 

Let $E:=E_{0}+E_{1}+\mathellipsis +E_{r}$. Then $K=EK\left[ + \right]\left( 
I-E \right)K$. Let us introduce $E_{r+1}:=I-E$, $K_{r+1}:=E_{r+1}K$, $\Gamma 
_{r+1}=E_{r+1}\Gamma $. Subspace $K_{r+1}$ is invariant with respect to $A$. 
From the construction of the Pontryagin space $K_{r+1}$ we conclude that all 
chains of $A$ at $\mathrm{\alpha }$ that are contained in $K_{r+1}$ are 
degenerate.

Using the above notation we can summarize the results in the following 
proposition.

\begin{proposition}\label{proposition4}Let $\alpha \in R$ be a generalized pole not of positive type of $Q\in N_{\kappa }(H)$ given by minimal representation (\ref{eq1}). Then 

\begin{enumerate}[(i)]%, (i), (ii),...
\item $\thinspace K=K_{0}\left[ + \right]K_{1}\left[ + \right]\mathellipsis \left[ + \right]K_{r}\left[ + \right]K_{r+1}\quad ,$ 

where $K_{i},\thinspace i=0,\thinspace 1,\thinspace \mathellipsis ,\thinspace 
r,\thinspace r+1$ are Pontryagin spaces of indexes $\kappa_{i}$, respectively; $\thinspace \kappa_{0}=0, \kappa 
=\sum\limits_{i=1}^{r+1} \kappa_{i} $. With $E_{i}$ denoting orthogonal projections onto $K_{i}$, operators $A_{i}:={AE}_{i}$ and $\Gamma 
_{i}:=E_{i}\Gamma $, are closely connected in ${K_{i}=E}_{i}K,\thinspace i=0,\thinspace 1,\thinspace 
2,\thinspace \mathellipsis ,\thinspace r+1$. For $i=1,\thinspace 2,\thinspace 
\mathellipsis ,\thinspace r$, operators $A_{i}$ have single eigenvalue $\alpha $ and respective non-degenerate Jordan chain ${\thinspace 
x}_{0}^{i},\thinspace \mathellipsis ,\thinspace x_{l_{i}-1}^{i}=\Gamma 
_{i}h_{i}$. Hilbert space $K_{0}$ consists of all positive eigenvectors of $A$ at $\alpha. $

\item For every $\thinspace i=1,\thinspace 2,\thinspace \mathellipsis ,\thinspace r$ there exist $h_{i}\in H$ such that subspace $K_{i}$ is a linear span of the Jordan vectors 
\[
x_{k}^{i}=\left( A-\alpha \right)^{l_{i}-1-k}\Gamma_{i}h_{i},
\quad
k=0,\thinspace \mathellipsis ,l_{i}-1.
\]
\item $\thinspace Q:=Q_{0}+Q_{1}+\mathellipsis +Q_{r}+Q_{r+1}\quad ,$ 
\end{enumerate}
\[
Q_{i}\left( z \right)=\Gamma_{i}^{+}\left( A_{i}-zE_{i} \right)^{-1}\Gamma 
_{i}\in N_{\kappa_{i}}\left( H \right),\thinspace i=1,\thinspace 
\thinspace \mathellipsis ,\thinspace r.
\]
\end{proposition}
Obviously, the decomposition obtained in the Proposition \ref{proposition4} is desirable 
and within $\alpha $.

\section{Inverse of the function $\mathbf{Q\left({z}\right)=}\mathbf{\Gamma}_{\mathbf{0}}^{\mathbf{+}}\left(\mathbf{A-z}\right)^{\mathbf{-1}}\mathbf{\Gamma }_{\mathbf{0}}$}\label{s4}
\begin{lemma}\label{lemma5} Let operators $\Gamma, \Gamma^{+}, J$ be as introduced in the Section \ref{s1}. Assume also that $\Gamma^{+}\Gamma $ is an boundedly invertible operator in the Hilbert space $\left( 
H,\left( .,. \right) \right)$. Then for operator 
\[
P:=\Gamma \left( \Gamma^{+}\Gamma \right)^{-1}\Gamma^{+}
\]
following statements hold:

\begin{enumerate}[(i)]%, (i), (ii),...
\item P is orthogonal projection in Pontryagin space $(K,\thinspace \left[ .,. \right]);$
\item Scalar product does not degenerate on $\Gamma \left( H \right)$ and therefore it does not degenerate on ${\Gamma \left( H \right)\thinspace }^{[\bot ]}=Ker\Gamma^{+}. $
\item $Ker\Gamma^{+}=\left( I-P \right)K$ 
\item Pontryagin space $K$ can be decomposed as a direct orthogonal sum of Pontryagin spaces i.e. 
\end{enumerate}
\begin{equation}
\label{eq19}
K=\left( I-P \right)K[+]PK.
\end{equation}
\end{lemma}
\textbf{Proof.} (i) Obviously $P^{2}=P$. 

According to well known properties of adjoint operators (see e.g. [JKL] p. 
34) it is easy to verify $\left[ \left( \Gamma^{+}\Gamma \right)^{-1} 
\right]^{\ast }= \quad \left( \Gamma^{+}\Gamma \right)^{-1}$ and then to verify 
$\left[ Px,y \right]=\left[ x,Py \right],\thinspace ie.\thinspace P^{[\ast 
]}=P$. This proves (i). 

(ii) If $\Gamma h\ne $ and $\left[ \Gamma h,\Gamma g \right]=0$, $\forall 
g\in H$, then $({\thinspace \Gamma }^{+}\Gamma hg)=0$, $\forall 
g\in H$. This means 

${\thinspace \Gamma }^{+}\Gamma h=0\to h=0\to \Gamma h=0.$ This is a 
contradiction that proves (ii).

(iii) It is sufficient to prove $Ker\Gamma^{+}=KerP$. 
\[
P:=\Gamma \left( \Gamma^{+}\Gamma \right)^{-1}\Gamma^{+}\to Ker\Gamma 
^{+}\subseteq KerP.
\]
Conversely, as $\Gamma^{+}\Gamma $ is boundedly invertible, 
\[
y\in KerP\to 0=\left[ \Gamma \left( \Gamma^{+}\Gamma 
\right)^{-1}\Gamma^{+}y,x \right]=\left( \left( \Gamma^{+}\Gamma 
\right)^{-1}\Gamma^{+}y,\Gamma^{+}x \right),\forall \Gamma^{+}x\in 
H.
\]
\[
R\left( \Gamma^{+} \right)=H
\quad
\to 
\quad
\left( \Gamma^{+}\Gamma \right)^{-1}\Gamma^{+}y=0\to \Gamma^{+}y=0\to 
y\in Ker\Gamma^{+}.
\]
(iv) Note that it holds $P\Gamma =\Gamma $ and $\Gamma^{+}P=\Gamma^{+}$. 
Now the statement (iv) follows directly from (iii) and (ii).
\\

If a function $Q$ is given by (\ref{eq7}) we define
\begin{equation}
\label{eq20}
P:=\Gamma_{0}\left( \Gamma_{0}^{+}\Gamma_{0} \right)^{-1}\Gamma 
_{0}^{+},
\end{equation}
\[
\tilde{A}:=\left( I-P \right)A\left( I-P \right):\left( I-P \right)K\to 
\left( I-P \right)K
\]
\[
{(\tilde{A}-z)}^{-1}:\left( I-P \right)K\to \left( I-P \right)K.
\]
Note 
\[
\left( I-P \right)\left( \tilde{A}-z \right)^{-1}\left( I-P \right)=\left( 
{\begin{array}{*{20}c}
\left( \tilde{A}-z \right)^{-1} & 0\\
0 & 0\\
\end{array} } \right)
\]
We prefer to use the notation on the left hand side, because it makes the 
following proofs shorter. 

\begin{theorem}\label{theorem2} Assume that function $Q\in N_{\kappa }(H)$ has the representation, 
\begin{equation}
\label{eq21}
Q(z)=\Gamma_{0}^{+}\left( A-z \right)^{-1}\Gamma_{0}
\end{equation}
where A is a self-adjoint bounded operator in a Pontryagin space K and $\Gamma_{0}^{+}\Gamma_{0}$ is boundedly invertible. Then the inverse function 
\[
\hat{Q}\left( z \right):=-{Q(z)}^{-1},\thinspace 
\]
has the following representation
\begin{equation}
\label{eq22}
\hat{Q}\left( z \right)=\left( \Gamma_{0}^{+}\Gamma_{0} \right)^{-1}\Gamma 
_{0}^{+}\left\{ A(I-P)\left( \tilde{A}-z \right)^{-1}(I-P)A-\left( A-z 
\right) \right\}\Gamma_{0}\left( \Gamma_{0}^{+}\Gamma_{0} \right)^{-1}.
\end{equation}
Note that we did \textbf{not} assume here that $Q$ satisfies minimality 
condition. 
\end{theorem}
\textbf{Proof. }For projection $P$ introduced by (\ref{eq20}), according to Lemma 
\ref{lemma5} (iv), we have the following decomposition 
\begin{equation}
\label{eq23}
A-zI=\left( {\begin{array}{*{20}c}
\left( I-P \right)(A-zI)(I-P) & \left( I-P \right)AP\\
PA(I-P) & P\left( A-zI \right)P\\
\end{array} } \right)\quad .
\end{equation}

For $z\in \rho (A)$ let us denote 
\begin{equation}
\label{eq24}
\left( {\begin{array}{*{20}c}
X & Y\\
Z & W\\
\end{array} } \right):={(A-z)}^{-1}\thinspace .
\end{equation}
By solving operator equations derived from the identity 
\[
\left( {\begin{array}{*{20}c}
\tilde{A}-z(I-P) & \left( I-P \right)AP\\
PA(I-P) & P\left( A-zI \right)P\\
\end{array} } \right)\left( {\begin{array}{*{20}c}
X & Y\\
Z & W\\
\end{array} } \right)=\left( {\begin{array}{*{20}c}
I-P & 0\\
0 & P\\
\end{array} } \right)\quad ,
\]
we get 
\begin{equation}
\label{eq25}
W=\left\{ P\left( A-zI \right)P-PA(I-P)\left( \tilde{A}-z 
\right)^{-1}(I-P)AP \right\}^{-1}.
\end{equation}
We need not to find operators $X$, $Y$ and $Z$. We only need to understand 
what their domains and ranges are. Then from 
\[
Ker(\Gamma_{0}^{+})=R(I-P),
\quad
R\left( \Gamma_{0} \right)=R(P)
\]
we easily see 
\[
\Gamma_{0}^{+}X\Gamma_{0}=\Gamma_{0}^{+}Y\Gamma_{0}=\Gamma 
_{0}^{+}Z\Gamma_{0}=0.
\]
By substituting (\ref{eq24}) into (\ref{eq21}) and using (\ref{eq25}) we get
\begin{equation}
\label{eq26}
Q\left( z \right)=\Gamma_{0}^{+}\left( {\begin{array}{*{20}c}
0 & 0\\
0 & W\\
\end{array} } \right)\Gamma_{0}=\Gamma_{0}^{+}\left\{ P\left( A-zI 
\right)P-PA(I-P)\left( \tilde{A}-z \right)^{-1}(I-P)AP \right\}^{-1}\Gamma 
_{0}.
\end{equation}
Then, by substituting (\ref{eq26}) and (\ref{eq22}) into the following product and using 
definition of $P$ we verify 
\[
Q\left( z \right)\hat{Q}\left( z \right)=
\]
\[
=\Gamma_{0}^{+}\left\{ P\left( A-zI \right)P-PA(I-P)\left( \tilde{A}-z 
\right)^{-1}(I-P)AP \right\}^{-1}\Gamma_{0}\times
\]
\[
\left( \Gamma_{0}^{+}\Gamma_{0} \right)^{-1}\Gamma_{0}^{+}\left\{ 
A(I-P)\left( \tilde{A}-z \right)^{-1}(I-P)A-\left( A-z \right) 
\right\}\Gamma_{0}\left( \Gamma_{0}^{+}\Gamma_{0} \right)^{-1}=
\]
\[
=\Gamma_{0}^{+}\left\{ P\left( A-zI \right)P-PA(I-P)\left( \tilde{A}-z 
\right)^{-1}(I-P)AP \right\}^{-1}\times 
\]
\[
\left\{ PA(I-P)\left( \tilde{A}-z \right)^{-1}(I-P)AP-P\left( A-zI \right)P 
\right\}\Gamma_{0}\left( \Gamma_{0}^{+}\Gamma_{0} \right)^{-1}
\]
\[
=\Gamma_{0}^{+}\left( -P \right)\Gamma_{0}\left( \Gamma_{0}^{+}\Gamma 
_{0} \right)^{-1}=-I\quad .
\]
We will later use representation (\ref{eq22}) to prove Theorem \ref{theorem3}, a result about desirable decomposition. Let us first give some consequences of the 
representation (\ref{eq22}). 

\begin{corollary}\label{corollary1}Let $Q\left( z \right),\thinspace \hat{Q}\left( z 
\right),\Gamma_{0},\thinspace \Gamma_{0}^{+}$ be the same as in Theorem \ref{theorem2}, then it holds 
\begin{equation}
\label{eq27}
\hat{Q}\left( z \right)\Gamma_{0}^{+}=\left( \Gamma_{0}^{+}\Gamma_{0} 
\right)^{-1}\Gamma_{0}^{+}\left\{ -I+A(I-P)\left( \tilde{A}-z 
\right)^{-1}(I-P) \right\}(A-zI).
\end{equation}
\end{corollary}
\textbf{Proof. }In the following derivations we will frequently use $\Gamma 
_{0}^{+}P=\Gamma_{0}^{+},\thinspace P\Gamma_{0}=\Gamma_{0}$. 
\[
\hat{Q}\left( z \right)\Gamma_{0}^{+}
\quad
=\left( \Gamma_{0}^{+}\Gamma_{0} \right)^{-1}\Gamma_{0}^{+}\left\{ 
A(I-P)\left( \tilde{A}-z \right)^{-1}(I-P)A-\left( A-zI \right) 
\right\}\Gamma_{0}\left( \Gamma_{0}^{+}\Gamma_{0} \right)^{-1}\Gamma 
_{0}^{+}
\]
\[
=\left( \Gamma_{0}^{+}\Gamma_{0} \right)^{-1}\Gamma_{0}^{+}\left\{ 
A(I-P)\left( \tilde{A}-z \right)^{-1}(I-P)(A-zI)P-\left( A-zI \right)P 
\right\}
\]
\[
=\left( \Gamma_{0}^{+}\Gamma_{0} \right)^{-1}\Gamma_{0}^{+}\times
\]
\[
\left\{ 
A\left( I-P \right)\left( \tilde{A}-z \right)^{-1}\left( I-P \right)\left( 
A-zI \right)\left( P-I \right)
+A(I-P)\left( \tilde{A}-z \right)^{-1}(I-P)\left( A-zI \right)-\left( A-zI \right)P \right\}
\]
\[
=\left( \Gamma_{0}^{+}\Gamma_{0} \right)^{-1}\Gamma_{0}^{+}\left\{ 
-A\left( I-P \right)+A(I-P)\left( \tilde{A}-z \right)^{-1}(I-P)\left( A-zI 
\right)-\left( A-zI \right)P \right\}
\]
\[
=\left( \Gamma_{0}^{+}\Gamma_{0} \right)^{-1}\Gamma_{0}^{+}\left\{ 
-\left( A-zI \right)+A(I-P)\left( \tilde{A}-z \right)^{-1}(I-P)\left( A-zI 
\right) \right\}
\]
\[
=\left( \Gamma_{0}^{+}\Gamma_{0} \right)^{-1}\Gamma_{0}^{+}\left\{ 
-I+A(I-P)\left( \tilde{A}-z \right)^{-1}(I-P) \right\}(A-zI).
\]
\begin{corollary}\label{corollary2} Let again $Q\left( z \right),\thinspace \hat{Q}\left( z 
\right),\Gamma_{0},\thinspace \Gamma_{0}^{+}$ be the same as in Theorem \ref{theorem2} Then the inverse $\hat{Q}$ has representation (\ref{eq4}), i.e.
\[
\hat{Q}\left( z \right)={-Q\left( \bar{z_{0}} 
\right)}^{-1}+(z-\bar{z_{0}})\hat{\Gamma }^{+}\left( \left( I+(z-z_{0} 
\right)\left( \hat{A}-z \right)^{-1} \right)\hat{\Gamma },
\]
where $\hat{A}$ is a self-adjoint linear \textbf{relation} with critical eigenvalue at $\infty $, i.e. it holds 
\begin{equation}
\label{eq28}
\hat{A}\left( 0 \right)=R\left( P \right)=R(\Gamma_{0}).
\end{equation}
\end{corollary}
\textbf{Proof.} The function $Q\in N_{\kappa }(H)$ that admits 
representation (\ref{eq21}) is a special case of the function that admits 
representation (\ref{eq1}). Let $z_{0}\in \rho (A)$ and let us introduce 
$\Gamma $ by 
\begin{equation}
\label{eq29}
\Gamma :=\left( A-z_{0} \right)^{-1}\Gamma_{0}\quad .
\end{equation}
Then from (\ref{eq21}) it easily follows (\ref{eq1}) where ${Q(z_{0})}^{\ast }=\Gamma 
_{0}^{+}\left( A-\bar{z_{0}} \right)\Gamma_{0}$. 

According to Lemma \ref{lemma2}, the inverse $\hat{Q}$ admits representation (\ref{eq4}). 
Then for $z=z_{0}$, from (\ref{eq5}) we get $\Gamma_{z_{0}}=\Gamma $, and from 
(\ref{eq6}) we get
\begin{equation}
\label{eq30}
\left( \hat{A}-z_{0} \right)^{-1}=\left( A-z_{0} \right)^{-1}-\Gamma 
_{z_{0}}{Q\left( z_{0} \right)}^{-1}\Gamma_{\bar{z_{0}}}^{+}
\end{equation}
From (\ref{eq29}), it follows

$\Gamma_{z_{0}}=\left( A-z_{0} \right)^{-1}\Gamma_{0}$ and $\Gamma 
_{\bar{z_{0}}}^{+}=\Gamma_{0}^{+}\left( A-z_{0} \right)^{-1}$.

Substituting this into (\ref{eq30}) gives
\[
\left( \hat{A}-z_{0} \right)^{-1}=\left( A-z_{0} \right)^{-1}-\left( A-z_{0} 
\right)^{-1}\Gamma_{0}{Q\left( z_{0} \right)}^{-1}\Gamma_{0}^{+}\left( 
A-z_{0} \right)^{-1}
\]
\[
=\left( A-z_{0} \right)^{-1}\left( I-\Gamma_{0}{Q\left( z_{0} 
\right)}^{-1}\Gamma_{0}^{+}\left( A-z_{0} \right)^{-1} \right)\quad .
\]
According to the Corollary \ref{corollary1} we get
\[
\left( \hat{A}-z_{0} \right)^{-1}=\left( A-z_{0} \right)^{-1}\left( 
I+P\left( -I+A(I-P)\left( \tilde{A}-z_{0} \right)^{-1}(I-P) \right) \right)
\]
\[
=\left( A-z_{0} \right)^{-1}\left( I-P+PA{\left( I-P \right)\left( 
\tilde{A}-z_{0} \right)}^{-1}(I-P) \right)
\]
\[
=\left( A-z_{0} \right)^{-1}\left( I+PA{\left( I-P \right)\left( 
\tilde{A}-z_{0} \right)}^{-1} \right)(I-P)
\]
From this we conclude $Ker\left( \hat{A}-z_{0} \right)^{-1}\supseteq R(P)$ 
and, therefore $\hat{A}\left( 0 \right)\supseteq R(\Gamma_{0})$. 

In order to prove,$\thinspace Ker\left( \hat{A}-z_{0} \right)^{-1}\subseteq 
R(\Gamma_{0})$, assume to the contrary that there exist 

$0\ne (I-P)y\in Ker\left( \hat{A}-z_{0} \right)^{-1}.$ Then from 
\begin{equation}
\label{eq31}
\left( \hat{A}-z_{0} \right)^{-1}=\left( A-z_{0} \right)^{-1}\left( I+PA{\left( I-P \right)\left( 
\tilde{A}-z_{0} \right)}^{-1} \right)(I-P)
\end{equation}
and from the fact that $z_{0}$ is a regular point of the operator $A$ we get
\[
\left( I+PA{\left( I-P \right)\left( \tilde{A}-z_{0} \right)}^{-1} 
\right)\left( I-P \right)y=0
\]
Let us define: ${\left( I-P \right)g:=\left( \tilde{A}-z_{0} 
\right)}^{-1}\left( I-P \right)y\ne0 $ 

Then we have 
\[
\left( \tilde{A}-z_{0} \right)\left( I-P \right)g+PA\left( I-P \right)g=0
\]
Having in mind $\tilde{A}=\left( I-P \right)A(I-P)$ we get 
\[
\left( \left( I-P \right)A-z_{0}(I-P)+PA \right)\left( I-P \right)g=0
\]
\[
\left( A-z_{0}\right)\left( I-P \right)g=0
\]
This means that $\left( I-P \right)g\ne0 $ is an eigenvector of $A$ in the 
eigenvalue $z_{0}$, which is a contradiction. Therefore, $Ker\left( 
\hat{A}-z_{0} \right)^{-1}=R(\Gamma_{0})$, which proves the statement.

\begin{remark}\label{remark3} One consequence of the Corollary \ref{corollary2} is that the 
inverse $\hat{Q}$ of $Q$ cannot have operator representation of the form (\ref{eq1}); $\hat{A}$ has to be a linear relation. Therefore, operator representation (\ref{eq22}) is 
essentially a new representation. 
\end{remark}
\begin{remark}\label{remark4}If $Q$ admits representation (\ref{eq21}) then a pole 
cancellation function can be conveniently defined by 
\begin{equation}
\label{eq32}
\eta \left( z \right):={Q\left( z \right)}^{-1}\Gamma_{0}^{+}\left\{ 
x_{0}+\left( z-\alpha \right)x_{1}+\mathellipsis +{(z-\alpha )}^{l-1}x_{l-1} 
\right\}\quad ,
\end{equation}
where $\left\{ x_{0},\thinspace \thinspace x_{1},\thinspace \mathellipsis 
,x_{l-1} \right\}$ is a Jordan chain of $\thinspace A$ at $\alpha $.

That is how the expression for $\hat{Q}\left( z \right)\Gamma_{0}^{+}$, 
proven in Corollary \ref{corollary1}, comes into play in the study of pole cancellation 
functions. 
\end{remark}
Pole cancelation functions of the form (\ref{eq32}) were constructed in \cite{B} for 
the functions $Q\in N_{\kappa }^{nxn}$ and were used there to characterize regular poles including their multiplicities. Existence of generalized poles was characterized in \cite{BL}, without characterization of their multiplicities. Much later, in \cite{BLu}, the functions of the form (\ref{eq32}) were used to characterize generalized poles of the function $Q\in N_{\kappa }(H)$, including their multiplicities. 

Note, if a Jordan chain of $\thinspace A$ at $\alpha $ of length $l$ 
saisfies $x_{l-1}=\Gamma_{0}h$, then, according to Corollary \ref{corollary1} the pole 
cancelation function (\ref{eq32}) has a very simple form $\eta \left( z 
\right)=-{(z-\alpha )}^{l}h$. 

\section{A desirable decomposition of the function $\hat{Q}$}\label{s5}
\begin{theorem}\label{theorem3} Let $Q(z)=\Gamma_{0}^{+}\left( A-z \right)^{-1}\Gamma 
_{0}\mathrm{\thinspace }\in N_{\kappa }(H)$, where A is a bounded self-adjoint operator and let $\alpha \in R$ be a generalized pole of $Q$. Assume that the derivative $Q^{'}\left( 
\infty \right):=\lim_{z\to \infty}{zQ(z)}=-\Gamma_{0}^{+}\Gamma_{0}$ is boundedly invertible operator. Then the inverse $\hat{Q}\left( z 
\right):=-{Q\left( z \right)}^{-1}$ has a desirable decomposition 
\[
\hat{Q}(z)=\hat{Q}_{1}(z)+\hat{Q}_{2}(z);
\quad
\hat{Q}_{1}\in N_{\hat{\kappa }_{1}}(H),
\quad
\hat{Q}_{2}\in N_{\hat{\kappa }_{2}}\left( H \right);
\quad
\hat{\kappa }_{1}+\hat{\kappa }_{2}=\kappa .
\]
That decomposition has the following properties:

\begin{enumerate} [(i)]%, (i), (ii), ...
\item Function $\hat{Q}_{1}$ may have a generalized zero at $\thinspace \alpha $ and cannot have any generalized pole in C.
\item The negative index $\hat{\kappa }_{1}$ is equal to the number of negative eigenvalues of the self-adjoint operator $\left( \Gamma_{0}^{+}\Gamma_{0} \right)^{-1}$ in the Hilbert space $H. $
\item If $Q$ has a generalized zero at $\alpha $ then the function $\hat{Q}_{2}$ has a generalized pole at $\alpha . $
\item The function $\hat{Q}_{2}$ has minimal representation $\hat{Q}_{2}\left( z \right):=\tilde{\Gamma }^{+}\left( \tilde{A}-z \right)^{-1}\tilde{\Gamma }, $
\end{enumerate}
where $\tilde{\Gamma }=(I-P)A\Gamma_{0}\left( \Gamma_{0}^{+}\Gamma_{0} 
\right)^{-1} , \tilde{A}=\left( I-P \right)A\left( I-P \right).$
\end{theorem}
\textbf{Proof.} From (\ref{eq22}), it follows $\hat{Q}=\hat{Q}_{1}+\hat{Q}_{2}$ 
where
\[
\hat{Q}_{1}\left( z \right)=-\left( \Gamma_{0}^{+}\Gamma_{0} 
\right)^{-1}\Gamma_{0}^{+}\left( A-z \right)\Gamma_{0}\left( \Gamma 
_{0}^{+}\Gamma_{0} \right)^{-1}
\]
\[
\hat{Q}_{2}\left( z \right):=\left( \Gamma_{0}^{+}\Gamma_{0} 
\right)^{-1}\Gamma_{0}^{+}A(I-P)\left( \tilde{A}-z \right)^{-1}(I-P)A\Gamma 
_{0}\left( \Gamma_{0}^{+}\Gamma_{0} \right)^{-1}
\]
Statements (i) and (iii) follow immediately. Let us prove the remaining 
statements. 

We know $\hat{Q}\in N_{\kappa }(H)$ and $\hat{\kappa }_{1}+\hat{\kappa 
}_{2}\ge \kappa $. Let us denote by $\kappa^{'}$ and $\kappa^{''}$ 
negative indexes of $PK$ and $\left( I-P \right)K$, respectively. Then, 
according to (\ref{eq19}) $\kappa^{'}+\kappa^{''}=\kappa $. 

For $f, g\in H$ we have
\[
\left( \frac{\hat{Q}_{1}\left( z \right)-{\hat{Q}_{1}\left( w \right)}^{\ast 
}}{z-\bar{w}}f,g \right)=\left( \left( \Gamma_{0}^{+}\Gamma_{0} 
\right)^{-1}f,g \right)
\]
As $\left( \Gamma_{0}^{+}\Gamma_{0} \right)^{-1}$ is bounded, hence 
defined on the whole $H$, we can consider $f=\Gamma_{0}^{+}\Gamma 
_{0}f_{0}$ and g$=\Gamma_{0}^{+}\Gamma_{0}g_{0}$, where $f_{0}$ and 
$g_{0}$ run through entire $H$ when $f, g$ run through entire $H$. Therefore 
\[
\left( \left( \Gamma_{0}^{+}\Gamma_{0} \right)^{-1}f,g \right)=\left( 
f_{0},\Gamma_{0}^{+}\Gamma_{0}g_{0} \right)=\left[ \Gamma_{0}f_{0},\Gamma 
_{0}g_{0} \right]
\]
As $R\left( \Gamma_{0} \right)=R(P)$ we conclude that $\hat{\kappa 
}_{1}=\kappa^{'}$. Therefore, $\kappa^{'}+\hat{\kappa }_{2}\ge \kappa 
=\kappa^{'}+\kappa^{''}$, and $\hat{\kappa }_{2}\ge \kappa^{''}$. 

If we introduce $\tilde{\Gamma }=(I-P)A\Gamma_{0}\left( \Gamma 
_{0}^{+}\Gamma_{0} \right)^{-1}$, then $\hat{Q}_{2}\left( z 
\right):=\tilde{\Gamma }^{+}\left( \tilde{A}-z \right)^{-1}\tilde{\Gamma }$, 
where $\tilde{A}$ is a self-adjoint operator in $\left( I-P \right)K$. That 
means $\hat{\kappa }_{2}\le \kappa^{''}$, i.e. $\hat{\kappa }_{2}=\kappa 
^{''}$. 

This proves the remaining statements of the theorem, including $\hat{\kappa 
}_{1}+\hat{\kappa }_{2}=\kappa $. 

\begin{corollary}\label{corollary3} Let again $Q\in N_{\kappa }\left( H \right)$ be as in Theorem \ref{theorem2} and let us consider representation (\ref{eq1}) of $\hat{Q}$. Let $(\hat{A},\thinspace 
\hat{\Gamma },\thinspace \hat{K})$ be a corresponding triplet. Then the linear relation $\hat{A}$ cannot be an operator. There exists an invariant subspace $\hat{K}_{1}\subseteq \hat{K}$ of the linear relation $\hat{A}$, where the negative index $\hat{\kappa 
}_{1}$ of $\hat{K}_{1}$ is equal to the number of negative eigenvalues of $\left( \Gamma_{0}^{+}\Gamma_{0} \right)^{-1}.$
\end{corollary}

\end{document}